\documentclass[12pt]{amsart}
\usepackage{amsmath,amsthm,epsfig,amssymb}
\usepackage{amsthm,amsfonts,amssymb,amscd}
\usepackage[all]{xy}
\usepackage{young}
\newtheorem{thm}[subsection]{Theorem}
\newtheorem{prop}[subsection]{Proposition}

\newtheorem{lem}[subsection]{Lemma}
\newtheorem{corol}[subsection]{Corollary}

\theoremstyle{definition}
\newtheorem{Def}[subsection]{Definition}
\newtheorem{Not}[subsection]{Notation}

\newtheorem{proposition-definition}[subsection]{Proposition-Definition}

\begin{normalsize} \end{normalsize}

\newcommand{\ZZ}{{\mathbb Z}}

\newcommand{\PP}{{\mathbb P}}
\newcommand{\NN}{{\mathbb N}}

\newcommand{\OOO}{{\mathcal O}}

\newcommand{\SSS}{{\mathcal S}}

\numberwithin{equation}{section}
\renewcommand\square{\frame{\phantom{{\large x}}}}

\author{F. Laytimi}

\address{F. L.: Math\'ematiques - b\^{a}t. M2, Universit\'e Lille 1,
F-59655 Villeneuve d'Ascq Cedex, France}
\email{fatima.laytimi@math.univ-lille1.fr}

\author{W. Nahm}

\address{W. N.: Dublin Institute for Advanced Studies,
10 Burlington Road, Dublin 4, Ireland}
\email{wnahm@stp.dias.ie}


\title{ A General Vanishing Theorem}

\begin{document}

\date{}

\begin{abstract}

Let $E$ be a vector bundle and $L$ be a  line bundle over a smooth
projective variety $X$.  In this article, we give a condition for
the vanishing of Dolbeault cohomology groups of the  form
$H^{p,q}(X,\SSS^{\alpha}E\otimes \wedge^{\beta} E\otimes L)$  when
$S^{\alpha+\beta}E \otimes L$ is ample. This condition is shown
to be invariant under the interchange of $p$ and $q$. The optimality
of this condition  is discussed for some parameter values.

\end{abstract}

\maketitle

\section{Introduction} \setcounter{page}{1}

Throughout this paper $X$  will denote
  a smooth projective variety of dimension $n$
over the field of complex numbers,
 $E$  a vector bundle of rank $e$, and  $ L$ a  line bundle  on $X. $

For any non-negative integers $\alpha ,\beta $ we denote by
$S^{\alpha} E, \ \wedge^{\beta}E$ the symmetric product and the
exterior product of $E.$ $H^{p,q}(X, S^{\alpha} E\otimes
\wedge^{\beta}E\otimes L )$ will denote the Dolbeault cohomology
group $$H^{q}(X, S^{\alpha} E\otimes \wedge^{\beta}E\otimes L\otimes
\Omega^p_X ),$$
  where $\Omega^p_X$ is the bundle of exterior differential forms of
   degree $p$ on $X.$

 We  start with some definitions.

\begin{Def}\label{def}
The function  $\delta:\NN\cup\{0\} \longrightarrow \NN$ is the one which
satisfies
$$\delta (x)=m \Longleftrightarrow {m\choose 2} \leq x < {{m+1}
\choose 2}. $$ 
The last two inequalities imply
$$\delta (x)=[\frac{\sqrt{8x+1}+1}{2}],
$$
 where the symbol [\ ] denotes
the integral part.\\

i.e., $\delta(0)=1, \ \delta(1)=\delta(2)=2,
\delta(3)=\delta(4)=\delta(5) =3, \\
\delta(6)=\delta(7)=\delta(8)=\delta(9)=4 , \ \ldots $

\end{Def}

\begin{thm} \label{main} Let $\alpha,\beta\in\NN$. If
 $S^{\alpha+\beta}E \otimes L$ is  ample , then
$$H^{p,q}(X, S^{\alpha}E\otimes \wedge^{\beta}E\otimes L)=0$$
$$\mbox {\rm{ for }} \ \ q+p-n>(r_0+\alpha)
(e+\alpha-\beta)-\alpha(\alpha+1),$$ $$\mbox {\rm{ where }}\ \
r_0=\rm{min} \{\beta, \delta (n-p ),\delta(n-q) \}.
$$
\end{thm}

\medskip

\begin{corol}\label{corol1} Let $\beta $ be a positive  integer.
If \ $S^{\beta}E \otimes L$ is  ample, then
$$H^{p,q}(X,\wedge^{\beta }E\otimes L)=0, $$
$$  \mbox{\quad  \rm {for}} \quad q+p-n> r_0(e-\beta). $$
$$\mbox {\rm{ where }}\ \
r_0=\rm{min} \{\beta, \delta (n-p ),\delta(n-q) \}.$$
\end{corol}

This Corollary  improve   the result of Manivel "theorem 1. p.91" in
\cite{M1}.

\medskip

\begin{corol}\label{corol2}
Assume $S^{\alpha}E \otimes L$ is ample. Then
$$H^{p,q}(X,S^{\alpha}E\otimes L)=0, $$
$$  \mbox{\quad  \rm {for}} \quad q+p-n> \alpha(e-1).$$

\end{corol}

This article is the final version of several attempts  \cite{N}, \cite{prep}.
The result of these  latest were  used by Chaput in \cite{Chap} and by
Laytimi-Nagaraj in \cite{Kyoto}.

In \cite{M3} Manivel studied the vanishing of Dolbeault cohomology of
a product of vector bundles tensored  with certain power of their determinant.
 The  presence of the latest allowed to deal with the problem by more direct method.

\section{The Schur Functor Version of the theorem}

Our main result is a consequence of a Schur functor version of the theorem, but
before giving this version, we need to recall some definitions and
results:

We start by  some preparation on partitions and Schur functors (for
a definition see  ~\cite {FH}).

A partition $u = (u_1,u_2,\dots ,u_r)$ is a sequence of non
increasing positive integers $u_i$. Its length is $r$ and its weight
is $|u|= \sum\limits_{i=1}^r u_i$. For $i>r$ we put  $u_i=0$. The
zero-partition is the one where all $u_i$ are zero.

For any partition $u$ the corresponding Schur functor is denoted by
$\SSS_u$.

Let $V$ be a vector space of dimension $d$. To each partition $u$
corresponds an irreducible $Gl(V)$-module $\SSS_u(V)$ which vanishes
iff $u_{d+1}>0$. For example, $\SSS_{(k)}V=S^k V$. By functoriality
the definition of Schur functors carries over to vector bundles  $E$
on $X$.

By abuse of language we say that $\SSS_u$ has a certain property, if
$u$ has this property. For example we will say  $S^k$ has weight $k$.

\begin{Def} \label{def2.2}
The Young diagram $Y(u)$ of a partition $u$ is given by
$$Y(u)=\{(i,j)\in\NN^2\mid j\leq u_i\}.$$
The transposed partition $\tilde u$ is defined by
$$Y(\tilde u)=\{(i,j)\in\NN^2\mid(j,i)\in Y(u)\}.$$
We use the notation $\wedge_u = \SSS_{\tilde u}.$
\end{Def}

\begin{Def}
The rank of a partition $u$, is $$rk(u)=max \{ \rho \mid
(\rho,\rho)\in Y(u)\}. $$ If $rk(u)=1$, then $u$ is called a hook.
\end{Def}

\begin{Not}\label{gamma} If $u$ is a hook with $u_1=\alpha+1$ and $|u|=k$, we write
$$\SSS_u=\Gamma^\alpha_k.$$

In particular, $\Gamma^0_k=\wedge^k$ and $\Gamma^{k-1}_k=S^k$.
\end{Not}

Recall that
$$S^{\alpha}E\otimes \wedge^\beta E=\Gamma^\alpha_{\alpha+\beta} E
\oplus\Gamma^{\alpha-1}_{\alpha+\beta} E.$$

\begin{Def}
For partitions $u,v$ of the same  weight, the dominance partial
ordering is defined by
$$u\succeq v,  \mbox { \ \  iff \ \ } \sum_{i=1}^j
u_i \ \geq \  \sum_{i=1}^j v_i  \mbox {\ \ for all \ } j.$$

This partial ordering can be extended to a pre-ordering of the set
of all non-zero partitions of arbitrary weight $u,v$ with
$|u|=n,|v|=m , $ by comparing as above  the partitions of the same
weight $m u$ and $n v, $ where the multiplication
 $$m u  = m\ (u_1, u_2, \ldots, u_r)=(m u_1,\ldots, m u_r) \
 \forall \  m   \in \NN.$$

\medskip

More precisely \ \ \ \ \ \ $ u\succeq v \mbox {\ \ \ \  iff \ \ \ \
} m u\succeq\ n v. $

\medskip

We write \ \ \ \ \ \ \ \ \ \ \ \   $u \simeq v {\mbox {\ \ \ \   iff
\ \ \ \   }}u \succeq v{\mbox {{\ } and {\ } }} v \preceq u.$

\end{Def}

When it is more convenient we will  write
 $\SSS_u\succeq \SSS_v$ instead of $u\succeq v$.
For example, $\wedge^r \succ  \wedge^{r+1}, $ and $S^{\alpha} \simeq
S^{1}$ for any ${\alpha}\in \NN $.

\begin{lem} {\rm {(Dominance Lemma)}} \rm{(\cite{Manuscripta} "theorem
3.7")} \label{Dom0}
\medskip

For any partition $u$ and $v.$ 

If \ $u\succeq v,$ \ then \ $\SSS_u E $ \ ample \ $ \Longrightarrow
\SSS_v E$  \   ample.\\

For example: 
If $\wedge^2 E$ is ample, then $\wedge^3 E$ is ample.
 
\end{lem}

Now we give the Schur presentation  of the main theorem  under which
 the main theorem will be shown. With the notation \ref{gamma} we have:

\begin{thm}\label{hook} Let $k \in \NN.$  If
 $S^kE \otimes L$ is  ample , then $$H^{p,q}(X, \Gamma^\alpha_k
E\otimes L )=0,$$ $$\mbox {\rm {for}}\ \  q+p-n>(r_0+\alpha)
(e-k+2\alpha)-\alpha(\alpha+1),$$
$$\mbox {\rm {where}}\ r_0=\rm{min} \{\beta, \delta (n-p
),\delta(n-q) \} . $$
\end{thm}
\medskip
\begin{prop}
Theorem \ref{hook} is equivalent to Theorem \ref{main}
\end{prop}
{\it Proof:}
 Since $${\mathcal
S}^{\alpha}E\otimes \wedge^{k-\alpha}E= \Gamma^{\alpha}_k E\oplus
\Gamma^{\alpha-1}_k E, $$ we have only to show that  for
$1\leq\alpha\leq k-1$ the conditions of Theorem  \ref{main} imply the
vanishing of $H^{p,q}(X, \Gamma^{\alpha-1}_k E)$, but this is clear
since the function $(r_0+\alpha) (e-k+2\alpha)-\alpha(\alpha+1)$
 is  increasing in  $\alpha$.

$\hfill{\square}$

\section{Some  Technical Lemmas}
We start with some proprieties of the function $\delta$ defined in
\ref{def}.
\begin{lem}\label{delta}
For $\mu\in\NN,  x\in \NN $ such that\\ $(x+\mu \delta(x),x-\mu
\delta(x)) \in \NN\times \NN,$ we have
\begin{align*}
& 1)\ \ \delta(x+\delta(x))=\delta(x)+1 \\
& 2)\ \ \delta(x+\mu \delta(x))\leq  \delta(x)+\mu \\
& 3)\ \ \delta(x-\mu \delta(x))\leq  \delta(x)-\mu .
\end{align*}
\end{lem}

{\it Proof:} The first assertion  and the case  $\mu=1$ in 2) and 3)  are
obvious.\\
 For both remaining  assertions  we use
induction on $\mu.$

For 2) \\  $\delta(x+\mu
\delta(x))=\delta(x+\delta(x)+(\mu-1) \delta(x)), $ since \\
$\delta(x)\leq \delta(x+\delta(x))= \delta(x)+1, $ we have\\
$\delta(x+\delta(x)+(\mu-1) \delta(x))\leq
\delta(x+\delta(x)+(\mu-1) \delta(x+\delta(x)).$ \\

Now  induction hypothesis gives \\ $\delta(x+\delta(x)+(\mu-1)
\delta(x+\delta(x))\leq \delta(x+\delta(x))+\mu-1=\delta(x)+\mu.$\\

For 3) \\
$\delta(x-\mu \delta(x))=\delta(x-\delta(x)-(\mu-1) \delta(x)), $\\
$\delta(x-\delta(x)-(\mu-1) \delta(x))\leq
\delta(x-\delta(x)-(\mu-1) \delta(x-\delta(x)). $\\
Induction hypothesis gives \\
 $\delta(x-\delta(x)+(\mu-1)
\delta(x-\delta(x))\leq \delta(x-\delta(x))-(\mu-1). $\\
 Now since it is true for $\mu=1, $ we get  $\delta(x-\delta(x))-(\mu-1)\leq
\delta (x) -\mu.$ $\hfill{\square}$

\begin{Def}\label{orderfi} Let $\phi:\NN\times \NN \rightarrow
\NN\times \NN \times \NN $ the following injection
$\phi(x,\alpha)=(\phi_1(x,\alpha), \phi_2(x,\alpha),\phi_3(x,\alpha
)), $  where
\begin{align*}
&\phi_1(x,\alpha)=\delta(x)+\alpha\\
&\phi_2(x,\alpha)=\ x-{{\delta(x)}\choose
{2}}\\
&\phi_3(x,\alpha)=\alpha
\end{align*}
We define an order on the pairs $(x, \alpha) \in \NN\times \NN$ by
the lexicographic order on $\NN\times \NN \times \NN $ induced by
$\phi,  $ we denote this order by
$$(x', \alpha')\leq_{\phi} (x, \alpha)$$
The set $\NN\times \NN$ endowed with the above order will be
denoted:
\begin{equation}\label{U}
\{\NN\times \NN,\ \  \leq_{\phi}\} :=\frak U
\end{equation}
\end{Def}

\begin{lem}\label{order}
For $\mu \in \ZZ-\{0 \}$ and  $(x+\mu \delta(x), \alpha-\mu) \in
\NN\times \NN, $ then
$$(x+\mu\delta(x), \alpha-\mu)\leq_{\phi}   (x,\alpha ).$$
where the order \ $\leq_{\phi}$\   is given in  Definition
\ref{orderfi}.
\end{lem}
{\it Proof:} By Lemma \ref{delta}\ \  $\phi_1(x+\mu\delta(x),
\alpha-\mu)\leq \alpha+\delta(x). $\\

If \ $\delta (x+\mu \delta(x))=\mu+\delta (x)$, then
$$\phi_2(x+\mu \delta(x), \alpha-\mu)=x-{{\delta(x)}\choose
{2}}-{{\mu}\choose {2}}\leq   x-{{\delta(x)}\choose {2}}.$$ If
${{\mu}\choose {2}}=0, $ which means $\mu=1, $ then $$\phi_3(x+\mu
\delta(x), \alpha-\mu)= \alpha-1< \alpha .$$ $\hfill{\square}$

We need to use these following results

\begin{lem}\label{large En}
Let  $E$ an ample vector  bundle and $G$ an arbitrary vector bundle
on a projective variety $X. $ Then for sufficiently large enough $n$
$S^ nE \otimes G$ is ample.
\end{lem}

\begin{lem}{\bf {Bloch-Gieseker}} \label{BG}  \cite{BG}
Let $L$ be a line bundle on a  projective variety $X$ and $d$ be a 
positive integer. Then there exist a projective variety $Y$, a
finite surjective morphism $f:Y \rightarrow  X,$ and a line bundle
$M$ on $Y, $ such that $f^*L\simeq M^d.$
\end{lem}
\begin{lem} \label{remains1}
Let $p,q,n,f_1,\ldots,f_r$ be fixed positive integers and

\noindent $\alpha^1,\ldots,\alpha^r$ be fixed non-zero partitions.
If $H^{p,q}(X,\otimes _{i=1}^r{\mathcal S}_{\alpha^i}F_i)=0$ for all
smooth  projective varieties  $X$ of dimension $n$ and all ample
vector bundles $F_1,\ldots,F_r$ of ranks $f_1,\ldots,f_r$ on $X$,
then this vanishing statement remains true if one of the $F_i$ is
ample and the others are nef.
\end{lem}
{\it Proof:} We can reorder the $F_i$ such that $F_1$ is ample.
Let $E=F_1$ and $\alpha=\alpha^1$. Let $N$ be a sufficiently large
number such that $S^N E\otimes\det E^*$ is ample  (for the
existence of such $N$ see Lemma \ref{large En}, and let
$a=\sum_{i=2}^m |\alpha^i|$. By Lemma ~\ref{BG} we can find a finite
surjective morphism $f:Y\rightarrow X, $ and a line bundle $M$ on
$Y, $ such that $f^*(\det E)=M^{Na}.$ Then $E_a=f^*E\otimes (M^*)^a$
is ample since $S^N E_a$ is. We have
$$ f^*({\mathcal S}_\alpha E\otimes_{i=2}^m {\mathcal S}_{\alpha^i}F_i)=
{\mathcal S}_\alpha E_a \otimes_{i=2}^m {\mathcal
S}_{\alpha^i}F'_i,$$ where $F'_i=M^{|\alpha|}\otimes f^*F_i$ for
$i=2,\ldots,m$. All $F'_i$ are ample. To finish the proof, we
use ``lemma 10 in \cite{M2} which says,
For any vector bundle ${\mathcal F}$ on $X$ and any finite
 surjective morphism
 $f:Y \rightarrow  X,$ the vanishing of $H^{p, q}(Y, f^*{\mathcal F})$
 implies the vanishing of $H^{p, q}(X, {\mathcal F}). $

$\hfill{\square}$

\begin{lem} \label{remove} Fix $n,p,q,k,\alpha\in\NN$ and $t\in\ZZ
.$ Assume that $$H^{p,q}(X,\Gamma^\alpha_k E)$$ vanishes for all
smooth projective varieties $X$ of dimension $n$ and all ample
vector bundles $E$ of rank \ $e=k+t$ on $X$. Let \ $\alpha<k'<k$.
Then $H^{p,q}(X,\Gamma^\alpha_{k'}E')$ vanishes for all ample
vector bundles $E'$ of rank $e'=k'+t$ on $X$.
\end{lem}
{\it Proof:} For given $E'$, put $E=E'\oplus L^{\oplus(k-k')}$,
where $L$ is any ample line bundle. Since $\Gamma
^\alpha_{k'}E'\otimes L^{k-k'}$ is a direct summand of  $\Gamma
^\alpha_k E$, we have
$$H^{p,q}(X,\Gamma^\alpha_{k'} E'\otimes L^{k-k'})=0$$
for ample vector bundle $E'$ of rank $e'$ and ample line
bundle $L$. By Lemma \ref{remains1} , this vanishing result remains
true, when $L$ is replaced by the trivial line bundle.
$\hfill{\Box}$
\begin{corol} \label{remain} Assume that there is an integer $k_0$ such that
$$H^{p,q}(X,\Gamma^{\alpha}_k E) = 0 \quad  \mbox{if } \quad  k>k_0,$$
for any projective smooth  variety $X$  of dimension $n$ and any
ample vector bundle $E$ of rank $e, $ under the  condition  \
$C(n,p,q,\alpha,e-k).$  Then under this same  condition the  vanishing
remains true for all $k$.
\end{corol}
The Bloch-Gieseker lemma can be used in other way to generalize
vanishing theorems. In particular one has
\begin{lem} \label{O_X} Fix $n,p,q,e\in\NN$ and partitions
$u,v$ of the same weight. Assume that $H^{p,q}(X,\SSS_u E)$ vanishes
for all projective varieties $X$ of dimension $n$ and all vector
bundles $E$ of rank $e$ for which $\SSS_v E$ is ample. Let $L$
be a line bundle and $F$ a vector bundle of rank $e$. Then
$H^{p,q}(X,\SSS_u F\otimes L)=0$, if $\SSS_v F\otimes L$ is
ample.
\end{lem}
{\it Proof:} Let's denote $|u|=|v| =d. $ By Lemma ~\ref{BG} we can
find a finite surjective morphism $f:Y\rightarrow X, $ and a line
bundle $M $ on $Y,$  such that $f^*L=M^d.$ Then
\begin{equation}\label{analogous}
\SSS_v (f^*F\otimes M)= f^*(\SSS_v F\otimes L)\ \ \mbox {\rm{is
ample. }}
\end{equation}
Due to the analogous equation \eqref{analogous} for $\SSS_u$ one has
by assumption
$$H^{p,q}(Y,f^*(\SSS_u F\otimes L))=0, $$ and the vanishing of
$H^{p,q}(X,\SSS_u F\otimes L)$ follows  by using
``lemma 10 in \cite{M2}.

 $\hfill{\Box}$

The lemma applies for example  if $\SSS_v F$ is nef and $L$ is
ample.

\begin{corol}\label{L} \rm {To generalize  vanishing of type
$H^{p,q}(X,\SSS_{u} F\otimes L), $ \ from \ $L=\mathcal O_X$ to
arbitrary $L$, it suffices to use Lemma \ref{O_X}.}
\end{corol}
We need to recall
\begin{lem}\rm{(\cite {layt} \label{nef} "lemma 1.3")}
Let $X$ be a projective variety,  $E, \ F $ be vector bundles on $X.
$  If  $E$ is ample and $F$ nef , then $E\otimes F$ is ample.
\end{lem}

\section{The Borel-Le Potier Spectral Sequence}

To prepare the proof, we need a  lemma and some properties of the
Borel-Le Potier spectral sequence, which has been made a standard
tool in the derivation of vanishing theorems \cite{Dem}.

Let $E$ be a vector bundle over a smooth projective  variety $X, $
$\rm{dim}(X)=n. $  Let $Y=G_r(E)$ be the corresponding Grassmann
bundle and $Q$ be the canonical quotient bundle over $Y$.

\begin{lem}\label{lem5.1}
Let $l, r$ be  positive integer and $k=lr, $ if \ $\wedge^r E $ is
ample. Then for $P+q>n+r(e-r)$
$$H^{P,q}(G_r(E), \det Q^l)=0.$$
\end{lem}

{\it Proof:} Since $\det Q= \OOO_{\PP (\wedge ^rE)}(1)|_{G_r(E)}.$
Thus $\Lambda^rE $ ample implies that $\det Q$ is ample. One
conclude by using  Nakano-Akizuki-Kodaira vanishing theorem
\cite{NAK}. 
$\hfill{\Box}$

\begin{Def} \label{def5.2} Let $\pi:Y\rightarrow X$ be a morphism of
projective manifolds, $P$ a positive integer and ${\mathcal F}$ a
vector bundle over $Y$. The Borel-Le Potier spectral sequence $^PE$
given by the data $\pi, P, {\mathcal F}$ is the spectral sequence
which abuts to $H^{P,q}(Y,{\mathcal F})$ , it is obtained from the
filtration on $\Omega^P_Y \otimes {\mathcal F}$ which is induced by
the filtration
$$F^p(\Omega^P_Y)=\pi^*\Omega^p_X \wedge \Omega^{P-p}_Y$$
on the bundle $\Omega^P_Y$ of exterior differential forms of degree
$P$.
\end{Def}

The graded bundle which corresponds to the filtration on
$\Omega^P_Y$ is given by
$$F^p(\Omega^P_Y)/F^{p+1}(\Omega^P_Y)=\pi^*\Omega^p_X\otimes \Omega^{P-p}_
{Y/X},$$ where $\Omega^{P-p}_{Y/X}$ is the bundle of relative
differential forms of degree $P-p$. Thus the $E_1$ terms of $^PE$
have the form
$$^PE_1^{p,q-p}= H^q(Y,\pi^*\Omega^p_X\otimes\Omega^{P-p}_{Y/X}\otimes
{\mathcal F}). $$

These $E_1$ terms can be calculated as limits groups of the Leray
spectral sequence associated to the projection $\pi, $
$$^{p,P}E^{q-j,j}_{2,L}=H^{p,q-j}(X, R^j \pi_*(\Omega^{P-p}_{Y/X}\otimes
\mathcal F ))$$

Now we  consider the Borel-Le Potier spectral sequence which abuts
to $H^{P,q}(G_r(E),\det Q^l)$.

\begin{prop} \label{prop5.3} Let $\pi:G_r(E)=Y\rightarrow X,$ the
 $E_1$ terms of the Borel-Le Potier spectral sequence  given by
$\pi,  P, \det Q^l$ have the form
$$^PE_1^{p,q-p}=
\bigoplus_{u \in \ \sigma(P-p,r)} H^q(G_r(E),\SSS_{u} Q^*\otimes\det
Q^{l}\otimes \wedge_u S \otimes \pi^*\Omega^p_X).$$ Here $S$ is the
tautological sub-bundle of $\pi^*E$ over $Y$ and $\sigma(p,r)$ is
the set of partitions of weight $p$ and length at most $r$.
\end{prop}

{\it Proof:} One has $\Omega_{Y/X}=Q^*\otimes S$. Thus
$$\Omega^{P-p}_{Y/X}=\bigoplus_{u \in \sigma (P-p, r)} {\SSS_u Q^*}\otimes
\wedge_u S. \ \ \ \ \ \ \ \  \ \ \ \ \ \ \ \  \ \ \ \ \ \ \ \ \ \ \
\ \ \ \ \  \ \ \ \ \ \ \ \  \ \ \ \ \ \ \ \  \hfill{\Box}$$

 Obviously Leray spectral sequence
degenerates at the $E_{2,L}$ level.

 Using the corollary 1. in  (\cite{M1} page 94) of Bott
formula,  Manivel computes the $E_1$ terms under some condition on
$P,  $ \rm {(\cite{M1}  Proposition 3. page 96)}. He states his
result under the supplementary condition $e\geq k, $ which is not
necessary for the calculation.

\begin{prop}\label{Manprop}\cite{M1}

Assume $P\geq n+(l-1){{r+1}\choose {2}}-l(r-1), $ and $k=lr.$  Let
\begin{align*}
&\alpha(p)=\frac {(l-1)(r+1)}{2}-\frac {P-p}{2}\\
&j(p)=(l-1) { {r}\choose {2}}-(r-1) \alpha (p).
\end{align*}

Then the $E_1$ terms of the spectral sequence have the form

\begin{displaymath}
 ^P E_1^{p,q}=\left\{ \begin{array}{lll}
  & H^{p,q-j(p)}(X,\Gamma
^{\alpha, k}E)& \textrm{for} \ \ (n-p,\alpha(p))\in \frak U \\
& 0  & \textrm {otherwise, }
\end{array}\right.
\end{displaymath}
where the set $\frak U$ is defined in \eqref{U}.
\end{prop}

Note that the connecting morphisms of Borel-Le Potier spectral
sequence
$$d_m: \ \ ^PE_m^{p,q-p}\  \longrightarrow \ \ ^PE_m^{p+m,q-p+1-m}$$
 all vanish, unless $m$ is a multiple of $r $ since
under $d_m$ the integer $\alpha $  goes to the integer $ \alpha
+\frac{m}{r}.$

\section{Proof of the main theorem}

\medskip

Before giving the proof of the main theorem, we will first explain
the case $r_0= \beta $ in the main theorem, which corresponds to
Corollary \ref {Nagoya2} bellow.

We need to  recall these results

\bigskip

\begin{thm}{\cite{Nagoya}} {\label{Nagoya}}
Let $E_i$ be vector bundles, with  $rank(E_i)= e_i, $  over a smooth
projective variety $X$ of dimension $n, $ and let $L$ be a line
bundle on $X.$ If \ \  $\otimes_{i=1}^m\Lambda ^{r_i} E_i\otimes L$
 is ample, then
$$H^{p,q}(X, \otimes_{i=1}^m\Lambda ^{r_i} E_i\otimes L)=0 \ \ \mbox{\rm
{for}} \ \ p+q-n> \sum_{i=1}^m r_i(e_i-r_i). $$
\end{thm}

\begin{corol} {\label{Nagoya2}}  Let $E$ be a vector
bundle of rank $e, $ and let $L$ be a line  bundle on a smooth
projective variety $X $ of dimension \  $n.$  If $S^{\alpha+\beta}E
\otimes L$ is ample, then

$$H^{p,q}(X, S^{\alpha} E\otimes \Lambda ^{\beta }E \otimes
L)=0  \mbox{\quad \rm {for}} \quad q+p-n> \alpha(e-1)+\beta(e-\beta
).$$
\end{corol}

{\it Proof:} We will apply the Theorem \ref{Nagoya} to the vector
bundle\\
$\underbrace{E\otimes E \cdots\otimes E}_ {\alpha {\mbox\ times}}
\otimes  \Lambda^{\beta} E\otimes L, $ which $S^{\alpha} E\otimes
\Lambda ^{\beta }E\otimes L$ is a direct summand of.

Let's first show  this equivalence of  ampleness
 \begin{equation}\label{sim}
   S^{\alpha}E\otimes F  \simeq
\underbrace {E\otimes E \cdots \otimes E }_{\alpha{\mbox\
times}}\otimes F 
 \end{equation}

for any vector bundles $F. $

Indeed: For the first  direction, Note that  $S^{\alpha}E\otimes L $
is direct summand of \ $\underbrace{E\otimes E \cdots\otimes E}_
{\alpha {\mbox\ times}}
\otimes F. $\\
For the second  direction, Littlewood-Richardson rules gives,
$$\underbrace{E\otimes E \cdots \otimes E}_ {\alpha {\mbox\ times}}=
S^{\alpha }E\oplus \sum_{|\lambda|=\alpha}S_{\lambda}E, $$ we have
clearly $\alpha  \succ \lambda$ in the dominance partial order. 
Use Remark \ref{Dom0}  to conclude.

Now by Littlewood-Richardson rules
 $$ S^{\alpha} E \otimes \Lambda^{\beta} E=\oplus \ \SSS_{\nu} E,
 \  \mbox{\rm{with}}\  |\nu|= \alpha +\beta, $$
satisfying  $\SSS_{\nu} \prec S^{\alpha+\beta}.$ Thus the
ampleness of $S^{\alpha+\beta}E \otimes L$ implies the ampleness
of $S^{\alpha} E\otimes \Lambda ^{\beta }E \otimes L$  by Remark
\ref{Dom0}. Use the equivalence of ampleness \eqref{sim} to conclude.

$\hfill{\Box}$

Due to Remark \ref{L} one can prove our main theorem without $L. $

\bigskip

 We prove Theorem \ref{main} by
induction on $(n-p,\alpha) \in \frak U, $ where the set $\frak {U} $
is given in  Definition \ref{orderfi}.

 Assume
that the result is true for all pairs $(p',\alpha') $ such that
$$(n-p',\alpha')\leq_{\phi} (n-p,\alpha), $$ with respect to the order
introduced  in Definition \ref{orderfi}.
\medskip

 Choose $r=\delta(n-p). $  Let $l$ be arbitrary if $n=p$,
otherwise let
   $l\geq \frac {r\alpha+n-p} {r-1} $. Choose $P$ such that
 $\alpha(p)= \alpha, $
and consider the Borel-Le Potier spectral sequence. Then for $k=lr$
$$^P E_1^{p,q+j(p)-p}=H^{p,q}(X,\Gamma ^ {\alpha,k} E). $$

 When $m$ is a
multiple of $r$, the morphisms $d_m$ connect
 $^PE_1^{p,q+j(p)-p}$ with  $^PE_1^{p',q'+j(p')-p'} $ where
 for the terms on the right of $^PE_1^{p,q+j(p)-p}$
\begin{equation}\label{r}
p'=p+\mu r,\ \ \  q'=q+\mu (r-1)+1,
\end{equation}
 and
\begin{equation}\label{l}
p'=p-\mu r,\ \ \  q'=q-\mu (r-1)-1
\end{equation}
for the  terms on the left. Here $\mu $ is any positive integer.

\begin{lem} For any  integers  $p'$ and $q'$  of the form \eqref{r} or 
\eqref{l},  $$^PE_1^{p',q'+j(p')-p'}=0, $$ when $q>Q(n-p, \alpha), $
where $$Q(n-p,
\alpha)=n-p+(\delta(n-p)+\alpha)(e-k+2\alpha)-\alpha(\alpha+1). $$
\end{lem}

{\it Proof:} The assertion is trivially true for \  $\alpha(p')<0 $
\ or \ $e-k+\alpha(p')<0, $ \  such that we may assume
$e-k+2\alpha(p') \geq 0. $

We need to prove that the assertion  $q>Q(n-p, \alpha)$ implies the
assertion  $q'>Q(n-p', \alpha'). $

The terms on the right of $^PE_1^{p,q+j(p)-p}$ have $\alpha'=
\alpha+ \mu, $ and  the parameters in  \eqref{r},
 a straight calculation yields

$$Q(n-p, \alpha)-Q(n-p', \alpha')+\mu( \delta (n-p)-1)+1=$$
\begin{equation}\label{Q1}
(e-k+2(\alpha+\mu)(\delta(n-p)-\mu-\delta(n-p-\mu\delta(n-p))+\mu^2+1.
\end{equation}

\medskip
The terms on the left of $^PE_1^{p,q+j(p)-p}$ have $\alpha'= \alpha-
\mu, $  and the parameters in  \eqref{l},  the calculation yields

$$Q(n-p, \alpha)-Q(n-p', \alpha')-\mu( \delta (n-p)-1)-1 =$$
\begin{equation}\label{Q2}
(e-k+2(\alpha-\mu)(\delta(n-p)+\mu-\delta(n-p+\mu\delta(n-p))+\mu^2-1.
\end{equation}
\medskip

By Lemma \ref{delta}  both   terms  of \eqref{Q1} and \eqref{Q2} are
non negative and positive if $\mu \neq 1. $ Thus $q
> Q(n-p, \alpha)$  implies $q' > Q(n-p', \alpha(p')). $

\medskip

By Lemma \ref{order}\  $(n-p', \alpha(p'))\leq_{\phi} (n-p, \alpha),
$ such that the groups $^PE_1^{p',q'+j(p')-p'}$  vanish by induction
hypothesis. Thus all co-bordant morphisms of $^P E_1^{p,q+j(p)-p}$
vanish. This implies that $^P E_1^{p,q+j(p)-p}$ is a sub-factor of
$H^{P,q+j(p)-p}(Y,F),  $ where $F=\det (Q)^l. $

Recall that $P= p+(l-1){ {r+1}\choose{2}} -\alpha r$ \ and \
$\rm{dim} Y=n+r(e-r). $ Thus the condition $q> Q(n-p,\alpha) $ is
equivalent to
$$P+q+j(p)-\rm{dim}Y> \alpha(e-k+\alpha).$$

When the right hand side is non-negative, $H^{P,q+j(p)-p}(Y,F)=0$ by
Nakano-Kodaira-Akizuki vanishing theorem. Thereby
$$H^{p,q}(X, \Gamma ^\alpha _k E )=0  \ \ \mbox{\rm{for}} \ \
q>Q(n-p,\alpha).$$
 Remember that this proof was under the condition
$k=rl$ see Proposition \ref{Manprop},  but this condition can be
removed by Corollary \ref{remain}.

To get $r_0=\delta(n-q)$ in our theorem,  we interchange the role of
$p$ and $q$ at every stage of the proof, in particular we use
$r=\delta(n-q).$

\section{Optimality}

\begin{prop}
Let $\it G=Gr_{(r,d)} $ be the Grassmannian of all co-dimensional
$r$ subspaces of a vector space $V$ of dimension $d=f+r. $ Let $Q$
be the universal sub-bundle of rank $r$ on $\it G,  $ $\dim\it
G=n=fr. $

Then,  for  $q= n-f$, $\alpha = f-1$
$$H^{q}(\it G, S^{\alpha}Q \otimes Q \otimes \det Q \otimes K_X) \neq
0$$
\end{prop}

{\it Proof:}   Since
 $S^{\alpha+1}Q $ is direct summand of
$S^{\alpha}Q \otimes Q, $ it's enough to show $$H^{q}(\it G,
S^{\alpha+1}Q \otimes \det Q \otimes K_X) \neq 0.  $$

For the universal sub-bundle $\it S$ on $\it G, $ we have $ K_{\it
G}=((\det Q)^*)^{\otimes d}= \det {\it S}^{\otimes d}. $

Thus since \   $\alpha = f-1$ $$H^{q}(\it G, S^{f}Q \otimes \det Q
\otimes K_X)= H^{q}(\it G, S^{f}Q \otimes \det {\it S}^{\otimes
{(d-1)}}. $$

Now by Bott formula (see corollary 1. page 94 of \cite{M1})
$$H^{q}(\it G, S^{f}Q \otimes \det {\it S}^{\otimes
{(d-1)}}=\delta _{q,i((a,b)-c(d))} \SSS _{\psi(a,b)}V, $$ where
$$a= ({f} ,\underbrace{0, \cdots,   0}_ {r-1 {\mbox\ times}}),  \
\ b=(\underbrace{d-1, \cdots, d-1 }_ {d-r {\mbox\ times}}). $$

For any sequence $v=(v_1, v_2, \ldots )$ $$i(v)= \rm {card} \{
(i,j)\ / \  i<j,\  v_i< v_j \}, $$ where
$$\psi(v)=(v-c(d))^{\geq}+c(d),$$
$$c(d) = (1,2, \ldots ,d), $$

and  $(v)^{\geq}$ \ is the partition
 obtained by ordering the  terms of $v$ in non increasing order.

$$(a,b)=(f ,\underbrace{0, \cdots,   0}_ {r-1 {\mbox\ times}}),
\underbrace{d-1, \cdots, d-1 }_ {d-r {\mbox\ times}}).$$

$$((a,b)-c(d))=(f-1  ,-2,-3,\ldots, -r,f-2,f-3,\ldots, 0,-1), $$ we
get $i((a,b)-c(d))=f(r-1)= n-f, $  and   $$\psi(a,b)=( \underbrace{
f,f,\ldots, f}_{{d}{\mbox\ times}}).$$

Thus $\SSS_{\psi((a,b)}V=(\det V)^{\otimes f}. $

Note that the  non-vanishing example of the above  proposition
happens for the limit condition
$$q+p-n=(r_0+\alpha) (e+\alpha-\beta)-\alpha (\alpha +1),$$
$$\mbox {\rm{ where }}\ \ r_0=\rm{min} \{\beta, \delta (n-p
),\delta(n-q) \}.
$$

 $\hfill{\Box}$

\end{document}